\documentclass[12pt]{article}
\usepackage{amsmath,comment}
\usepackage{amssymb}
\usepackage{amsthm}
\providecommand{\abs}[1]{\lvert#1\rvert}

\newtheorem{thm}{Theorem}

\newtheorem{lem}{Lemma}

\newenvironment{proofof}[1]{\noindent{\itshape
    #1. }}{\hfill$\square$\medskip}
\newcommand{\dist}{\mbox{\em dist\/}}
\newcommand{\ce}{{\cal E}}
\newcommand{\ch}{{\cal H}}


\begin{document}

\hfill 220907

\begin{center}
{\large\bf METRIC HYPERGRAPHS AND METRIC-LINE EQUIVALENCES}\\
\rule{0pt}{17pt} Va\v sek Chv\' atal and Ida Kantor
\end{center}

\begin{abstract}
In a metric space $M=(X,\dist)$, we say that {\em $v$ is between $u$ and $w$} if $\dist(u,w)=\dist(u,v)+\dist(v,w)$. Taking all triples $\{u,v,w\}$ such that $v$ is between $u$ and $w$, one can associate a 3-uniform hypergraph  
with each finite metric space $M$. An effort to solve some basic open questions regarding finite metric spaces has motivated an endeavor to better understand these associated hypergraphs. In answer to a question posed in \cite{BBC13}, we present an infinite family of hypergraphs that are non-metric, i.e., they don't arise from any metric space. 

Another basic structure associated with a metric space is a binary equivalence on the vertex set, where two pairs are in the same class if they induce the same line. An equivalence that comes from some metric space is a {\em metric-line equivalence}. We present an infinite family of so called {\em obstacles}, that is, binary equivalences that prevent an equivalence from being a metric-line equivalence.
\end{abstract}

\section{Introduction}
Given a metric space $(V,\dist)$, we follow \cite{ACH} in writing $[uvw]$ to signify that 
$u,v,w$ are pairwise distinct points of $V$ and $\dist(u,v)+\dist(v,w)=\dist(u,w)$. With $M$ standing for the metric space, it will be convenient to write
\[
\ce_M \;=\; \{\{x,y,z\}:\; [zxy] \text{ or } [xzy] \text{ or } [xyz].\}.
\]
Following~\cite{BBC13}, we say that a $3$-uniform hypergraph $(V,\ce)$ is {\em metric\/}
if there is a metric space $M$ such that $\ce= \ce_M$. All induced subhypergraphs of metric hypergraphs are metric, and so metric hypergraphs can be characterized as hypergraphs without certain induced
subhypergraphs, namely, the minimal non-metric ones. Section 3 of~\cite{BBC13} presents three minimal non-metric hypergraphs along with the following comment:
\begin{quote}
If there are only finitely many minimal non-metric hypergraphs, then metric hypergraphs can
be recognized in polynomial time. However, it is conceivable that there are
infinitely many minimal non-metric hypergraphs and it is not clear whether
metric hypergraphs can be recognized in polynomial time.
\end{quote}
Our first main result shows that there are indeed infinitely many minimal non-metric hypergraphs. To formulate it, we need the following definition: When $G$ is a graph with vertex set $V$ and edge set $E$, the {\em hypergraph based on $G$\/} is the $3$-uniform hypergraph with vertex set $V\cup\{x\}$, where $x\not\in V$, and hyperedge set  consisting of all three-point subsets of $V$ and all three-point sets $\{x,u,v\}$ such that $\{u,v\}\in E$. 
\begin{thm}\label{evenc}
For every even integer $n$ greater than four, the hypergraph based on the cycle $C_n$ is minimal non-metric.
\end{thm}

In 1943, Erd\H os~\cite{Erd43} proved that a set of $n$ points in the Euclidean plane determines at least $n$ distinct lines unless these $n$ points are collinear. In 2006, Chen and Chv\' atal~\cite{CC08} asked whether the same statement holds true more generally in all metric spaces $M$ with the line $L_M(xy)$ determined by two points $x$ and $y$ defined by 
\[
L_M(xy) \;=\; \{ x,y\}\,\cup\,\{z: \{x,y,z\}\in \ce_M\}.
\]
Early progress toward the conjecture that this generalization does hold true is surveyed in \cite{Chv18}; contributions too recent to be included there are \cite{BKR19}, \cite{ABMZ22}, \cite{Kan20}, \cite{RV21}, \cite{SS20}.

When $W$ is a set, 
we let $\binom{W}{2}$ denote the set of all $2$-point subsets of $W$. We say that an equivalence relation 
$\equiv$ on $\binom{W}{2}$ is a {\em metric-line equivalence\/} if there is a metric space $M$ on ground set $W$ such that
\begin{equation}\label{meq}
L_M(ab)=L_M(cd) \;\Leftrightarrow\; \{a,b\}\equiv\{c,d\}
\end{equation}
for every choice of two-point subsets $\{a,b\}$ and $\{c,d\}$ of $W$. 

How difficult is it to tell which equivalence relations on $\binom{W}{2}$ are metric-line equivalences and which of them are not?
Attempts at answering the Chen-Chv\' atal question could be only helped by an efficient algorithm for their recognition. Until now, all we have had here was a polynomial-time algorithm that, given an equivalence relation $\equiv$ on $\binom{W}{2}$, will in some cases certify that $\equiv$ is not a metric-line equivalence~\cite[Algorithms G and H]{Chv18}. We are going to offer another such certificate.

We say that an equivalence relation $\equiv$ on $\binom{V}{2}$ is an {\em obstacle\/} if no metric space $M$ on a superset $W$ of $V$ satisfies \eqref{meq} for every choice of two-point subsets $\{a,b\}$ and $\{c,d\}$ of $V$. 
It may not be obvious that there exist any obstacles at all; our second main result shows that there are infinitely many genuinely different  ones. To formulate this result, we need additional definitions again.

We say that an obstacle $\equiv$ on $\binom{V}{2}$ is {\em minimal\/} if there is no proper subset $U$ of $V$ such that the restriction of $\equiv$ on $\binom{U}{2}$ is an obstacle. Given a graph $G$ with vertex set $V$ and edge set $E$, we define the equivalence relation $\stackrel{G}{\equiv}$ on $\binom{V}{2}$
by 
\[
e \stackrel{G}{\equiv} f \;\Leftrightarrow\; (e\in E, f\in E)\; \text{ or }\;(e\not\in E, f\not\in E)
\]

\begin{thm}\label{obs}
If $G=C_n$ with $n$ an even integer greater than four, then the equivalence relation $\stackrel{G}{\equiv}$ is a minimal obstacle.
\end{thm}

\section{Metric and non-metric hypergraphs}

A corollary of the following lemma is one of the ingredients of our proof of Theorem~\ref{evenc}. 
\begin{lem}\label{oddc}
For every odd integer $n$ greater than one, the hypergraph based on $C_n$ is metric.
\end{lem}
\begin{proof} Writing $n=2s+1$, consider the metric space on the ground set $\{0,1,\ldots ,2s\}\cup\{x\}$ with 
$\dist(i,j)=j-i$ whenever $0\le i<j<n$ and
\[
\dist(x,k)=\begin{cases}
s& \text{if $k$ is even,}\\
s+1& \text{if $k$ is odd}.
\end{cases}
\]
\end{proof}
The next lemma comes from~\cite{RiRi} and has been generalized as~\cite[Lemma~3.2]{BBC13}.
\begin{lem}\label{riri}
Let $M$ be a metric space and let $V$ be a subset of its ground set such that 
$\abs{V}\ge 5$. If every three-point subset of $V$ belongs to $\ce_M$, then 
the elements of $V$ can be renamed as $0$, $1$, \ldots $n-1$ with $n=\abs{V}$ in such a way that 
\begin{equation}\label{linear}
0\le u<v<w<n \;\Rightarrow\; [uvw]
\end{equation}
\end{lem}
\begin{lem}\label{wh}
Given a metric space $M$ on a ground set $V\cup\{x\}$, where $x\not\in V$ and $V=\{0,1,\ldots n-1\}$, set 
\begin{eqnarray*}
D_1 &=& \{(j,\ell ):\;j<\ell \text{ and } [xj\ell ]\},\\
D_2 &=& \{(j,\ell ):\;j<\ell \text{ and } [jx\ell ]\},\\
D_3 &=& \{(j,\ell ):\;j<\ell \text{ and } [j\ell x]\}.
\end{eqnarray*}
If \eqref{linear} holds true, then
\begin{eqnarray}
(j,\ell)\in D_1,\; j<k<\ell\;\;&\Rightarrow&\; \;(j,k)\in D_1,\;\; (k,\ell)\in D_1,\label{xjkl}\\
(j,\ell)\in D_3,\; j<k<\ell\;\;\,&\Rightarrow&\; \,(j,k)\in D_3,\;\; (k,\ell)\in D_3,\label{jklx}\\
(j,\ell)\in D_2,\; i<j<\ell\;\;\;&\Rightarrow&\; \,\;(i,j)\in D_3,\;\; (i,\ell)\in D_2,\label{ijxl}\\
(j,\ell)\in D_2,\; j<\ell<m\;&\Rightarrow&\; (j,m)\in D_2,\;\; (\ell, m)\in D_1, \label{jxlm}\\
(i,j)\in D_1,\; (j,k)\in D_1\;&\Rightarrow&\; (i,k)\in D_1,\label{11}\\
(i,j)\in D_3,\; (j,k)\in D_3\;&\Rightarrow&\; (i,k)\in D_3,\label{33}\\
(i,k)\in D_2,\; (j,k)\in D_1\;&\Rightarrow&\; (i,j)\;\in D_2\;\;\text{\rm or } (j,i)\in D_2.\label{last}\\
(i,k)\in D_2,\; (i,j)\in D_3\;&\Rightarrow&\; (j,k)\in D_2\;\;\text{\rm or } (k,j)\in D_2,\label{first}
\end{eqnarray} 
\end{lem}
\begin{proof}
Implications \eqref{xjkl}\ ---\ \eqref{last} are special cases of the easily verifiable 
\begin{equation}\label{pm}
[abd],\,[bcd]\;\;\Rightarrow\;\;[abc],[acd],
\end{equation}
which has been pointed out first by Menger~\cite{M28}.
\end{proof}
Actually, the conclusion of \eqref{last} can be strengthened to $(i,j)\in D_2$ since $j<i$ would contradict \eqref{xjkl} with $(j,i,k)$ in place of $(j,k,\ell)$. Similarly, the conclusion of \eqref{first} can be strengthened to $(j,k)\in D_2$ since $k<j$ would contradict \eqref{jklx} with $(i,k,j)$ in place of $(j,k,\ell)$. These niceties are irrelevant to our purpose.
\begin{lem}\label{new}
Let $G$ be a graph with the vertex set $V=\{0,1,\dots,n-1\}$ and an edge set $E$. Let $x$ be a point outside $V$ and let $M$ be a metric space on the ground set $V\cup\{x\}$ such that, for every $u,v,w\in V$,
\begin{itemize}
	\item  $u<v<w \;\Rightarrow\; [uvw]$
	\item $\{u,v\}\in E \Leftrightarrow \{x,u,v\}\in \ce_M$.
\end{itemize} 

If $G$ contains no triangle, then
\[E\subseteq\{\{0,1\},\{1,2\},\ldots \{n-2,n-1\},\{n-1,0\}\}.\]
\end{lem}
\begin{proof}
Lemma~\ref{wh} guarantees that \eqref{xjkl}\ ---\ \eqref{last} are satisfied. Since a two-point subset $\{u,v\}$ of $V$ belongs to $E$ if and only if $\{x,u,v\}\in \ce_M$, we have 
\[E=\{\{u,v\}: (u,v)\in D_1\cup D_2\cup D_3\}\]
If $G$ contains no triangle, then \eqref{xjkl}, \eqref{jklx} show that every $(j,\ell)$ in $D_1\cup D_3$ has $\ell=j+1$ and \eqref{ijxl} and \eqref{jxlm}
 show that every $(j,\ell)$ in $D_2$ has $j=0$, $\ell=n-1$. Therefore
\[D_1\cup D_3 \subseteq \{(i,i+1): 0\le i< n-1\} \text{\ and \ }
D_2 \subseteq \{(0,n-1)\}.
\]\end{proof}

\begin{proofof}{Proof of Theorem~\ref{evenc}}
First, we shall deduce a contradiction from the assumption that the hypergraph $\ch$ based on a $C_n$ with $n$ even and greater than four is metric. For this purpose, let $V$ denote the vertex set of the $C_n$ and let $E$ denote its edge set. 
If $\ch$ is metric, then Lemma~\ref{riri} guarantees that (possibly after renaming the vertices) the hypothesis of Lemma~\ref{new} with $G=C_n$ is satisfied, and so 
\[E=\{\{0,1\},\{1,2\},\ldots \{n-2,n-1\},\{n-1,0\}\}.\]
By \eqref{first} with $i=0$, $j=1$, $k=n-1$, we have
\[
(0,1)\in D_1;
\] 
by \eqref{last} with $i=0$, $j=n-2$, $k=n-1$, we have
\[
(n-2,n-1)\in D_3;
\] 
by \eqref{11} and \eqref{33} with $j=i+1$, $k=i+2$, we have
\[
(i,i+1)\in D_1 \; \Leftrightarrow \; (i+1,i+2)\in D_3\;\;\;\text{ whenever $0\le i<n-2$.}
\]
Therefore the ordered pairs $(0,1)$, $(1,2)$, $(2,3)$, \ldots, $(n-2,n-1)$ alternate between $D_1$ and $D_3$, beginning with 
$(0,1)$ in $D_1$ and ending with $(n-2,n-1)$ in $D_3$. This contradicts the assumption that $n$ is even.

It remains to be proved that every proper induced subhypergraph $\ch_0$ of $\ch$ is metric, For this purpose,  note that $\ch_0$ is an induced subhypergraph of the complete $3$-uniform hypergraph on $n$ vertices, which is trivially metric, or an induced subhypergraph of the hypergraph based on the path of order $n-1$, which (being a subhypergraph of the hypergraph based on any larger cycle) is metric by Lemma~\ref{oddc}.
\end{proofof}

The lower bound on $n$ in Theorem~\ref{evenc} is essential: the hypergraph based on $C_4$ is metric. To see this, note that the metric space on the ground set $\{a,b,c,d,x\}$ with metric defined by the chart
\[
\begin{array}{r|c|c|c|c|c|}
  & a & b & c & d & x\\ \hline
a & 0 & 1 & 2 & 1 & 2 \\ \hline
b & 1 & 0 & 1 & 2 & 3 \\ \hline
c & 2 & 1 & 0 & 1 & 2 \\ \hline
d & 1 & 2 & 1 & 0 & 3 \\ \hline
x & 2 & 3 & 2 & 3 & 0 \\ \hline
\end{array}
\]
has
\[
[abc], [bcd], [cda], [dab],\; [xab], [xad], [xcb], [xcd],
\]
but none of $[xac]$, $[axc]$, $[acx]$, $[xbd]$, $[bxd]$, $[bdx]$.

The following lemma is used in the next section in our proof of Theorem~\ref{obs}:
\begin{lem}\label{house}
The hypergraph based on the complement $\overline{P_5}$ of the path of order 
five is non-metric.
\end{lem}
\begin{proof}
We shall deduce a contradiction from the assumption that the hypergraph $\ch$ 
based on $\overline{P_5}$ is metric. 
For this purpose, let $V$ denote the vertex set of our $\overline{P_5}$ and 
let $E$ denote its edge set.
If $\ch$ is metric, then there is a metric space $M$ on the ground set $V\cup\{x\}
$, where $x\not\in V$, such that all three-point subsets of $V$ belong to $\ce_M
$ and such that a two-point subset $\{u,v\}$ of $V$ belongs to $E$ if and 
only if $\{x,u,v\}\in \ce_M$. 
By Lemma~\ref{riri}, the elements of $V$ can be renamed as $0$, $1$, \ldots, $4
$ in such a way that 
\[
0\le u<v<w<5 \;\Rightarrow\; [uvw]
\]
and Lemma~\ref{wh} guarantees that \eqref{xjkl}\ ---\ \eqref{first} are 
satisfied. Since a two-point subset $\{u,v\}$ of $V$ belongs to $E$ if and 
only if $\{x,u,v\}\in \ce_M$, we have 
\[
E=\{\{u,v\}: (u,v)\in D_1\cup D_2\cup D_3\}.
\]
Since $\overline{P_5}$ contains only one triangle, \eqref{xjkl} and \eqref{jklx} show that every $(j,\ell)$ in $D_1\cup D_3$ has $\ell\le j+2$ and 
\eqref{ijxl}, \eqref{jxlm} show that every $(j,\ell)$ in $D_2$ has $\ell-j\ge 3$. Explicitly, we have 
\begin{eqnarray*}
D_1\cup D_3 &\subseteq& \{ (0,1),(1,2),(2,3),(3,4),(0,2),(1,3),(2
,4)\},\\
D_2 &\subseteq& \{(0,3),(0,4),(1,4)\}.
\end{eqnarray*}
Next, let $T$ denote the unique triangle in our $\overline{P_5}$ and let $P$ denote the three-edge path resulting when the three edges of $T$ are deleted.

Since no edge of $P$ extends to a triangle, implications \eqref{xjkl},\eqref{jklx} show that none of $\{0,2\},\{1,3\},\{2,4\}$ can belong to $P$ and implications \eqref{ijxl},\eqref{jxlm} show that neither of $\{0,3\},\{1,4\}$ can belong to $P$. 
Hence each of the three edges of $P$ must be one of $\{0,1\},\{1,2\},\{2,3\},\{3,4\},\{0,4\}$.

If the vertex set of $T$ is $\{0,1,3\}$, then $(1,3)\in D_1\cup D_3$ and by \eqref{xjkl} or \eqref{jklx}, $E$ contains also $\{1,2\}$ and $\{2,3\}$. However, $E$ contains only one triangle, so this won't happen. We rule out the possibility of $T$ being one of $\{0,2,3\}$, $\{0,2,4\}$, $\{1,2,4\}$, $\{1,3,4\}$ similarly. The remaining options are $\{0,1,2\}$, $\{1,2,3\}$, $\{2,3,4\}$, $\{0,1,4\}$, and $\{0,3,4\}$. 
Flip symmetry $i\leftrightarrow 4-i$ reduces these to the following three.\\ 
\phantom{xxx}{\sc Option 1}: the vertex set of $T$ is $\{0,1,2\}$,\\
\phantom{xxx}{\sc Option 2}: the vertex set of $T$ is $\{1,2,3\}$,\\
\phantom{xxx}{\sc Option 3}: the vertex set of $T$ is $\{0,3,4\}$.\\
We are going to eliminate these three options one by one.

{\sc Option 1}: $E=\{\{0,1\},\{1,2\},\{0,2\},\{2,3\},\{3,4\},\{0,4\}\}$.\\
\phantom{xx}$(0,4)\in D_2$ and $\{2,4\}\not\in E$ force $(0,2)\not\in D_3$ [and so $(0,2)\in D_1$] by \eqref{first},\\ 
\phantom{xx}$(0,2)\in D_1$ and $\{0,3\}\not\in E$ force $(2,3)\not\in D_1$ [and so $(2,3)\in D_3$] by \eqref{11},\\
\phantom{xx}$(0,4)\in D_2$ and $\{0,3\}\not\in E$ force $(3,4)\not\in D_1$ [and so $(3,4)\in D_3$] by \eqref{last},\\ 
\phantom{xx}$(2,3)\in D_3$ and $(3,4)\in D_3$ force $(2,4)\in D_3$ by \eqref{33}.\\
However, $(2,4)\in D_3$ is incompatible with $\{2,4\}\not\in E$.

{\sc Option 2}: $E=\{\{1,2\},\{2,3\},\{1,3\},\{0,1\},\{3,4\},\{0,4\}\}$.\\
\phantom{xx}$(0,4)\in D_2$ and $\{1,4\}\not\in E$ force $(0,1)\not\in D_3$ [and so $(0,1)\in D_1$] by \eqref{first},\\ 
\phantom{xx}$(0,1)\in D_1$ and $\{0,3\}\not\in E$ force $(1,3)\not\in D_1$ [and so $(1,3)\in D_3$] by \eqref{11},\\
\phantom{xx}$(0,4)\in D_2$ and $\{0,3\}\not\in E$ force $(3,4)\not\in D_1$ [and so $(3,4)\in D_3$] by \eqref{last},\\ 
\phantom{xx}$(1,3)\in D_3$ and $(3,4)\in D_3$ force $(1,4)\in D_3$ by \eqref{33}.\\
However, $(1,4)\in D_3$ is incompatible with $\{1,4\}\not\in E$.

{\sc Option 3}: $E=\{\{0,3\},\{0,4\},\{3,4\},\{0,1\},\{1,2\},\{2,3\}\}$.\\
\phantom{xx}$(0,3)\in D_2$ and $\{1,3\}\not\in E$ force $(0,1)\not\in D_3$ [and so $(0,1)\in D_1$] by \eqref{first},\\ 
\phantom{xx}$(0,1)\in D_1$ and $\{0,2\}\not\in E$ force $(1,2)\not\in D_1$ [and so $(1,2)\in D_3$] by \eqref{11},\\
\phantom{xx}$(0,3)\in D_2$ and $\{0,2\}\not\in E$ force $(2,3)\not\in D_1$ [and so $(2,3)\in D_3$] by \eqref{last},\\ 
\phantom{xx}$(1,2)\in D_3$ and $(2,3)\in D_3$ force $(1,3)\in D_3$ by \eqref{33}.\\
However, $(1,3)\in D_3$ is incompatible with $\{1,3\}\not\in E$. 
\end{proof}
By the way, the hypergraph based on $\overline{P_5}$ is minimal non-metric. To verify this, enumerate the vertices of $\overline{P_5}$ as $a,b,c,d,e$ in such a way that the edges of this $\overline{P_5}$ are
\[
\{a,b\},\; \{b,c\},\; \{c,d\},\; \{d,a\},\; \{e,b\},\; \{e,c\}.
\]
Now each of $\overline{P_5}-b$ and $\overline{P_5}-c$ is a $P_4$; by Lemma~\ref{oddc}, the hypergraph based on $P_4$ is metric. 
Next, $\overline{P_5}-e$ is a $C_4$; by the comment following our proof of Theorem~\ref{evenc}, the hypergraph based on $C_4$ is metric. 
Finally, $\overline{P_5}-a$ and $\overline{P_5}-d$ are isomorphic; to see that the hypergraph based on $\overline{P_5}-a$ is metric, note that the metric space on the ground set $\{e,b,c,d,x\}$ with metric defined by the chart
\[
\begin{array}{r|c|c|c|c|c|}
  & e & b & c & d & x\\ \hline
e & 0 & 1 & 2 & 3 & 2 \\ \hline
b & 1 & 0 & 1 & 2 & 3 \\ \hline
c & 2 & 1 & 0 & 1 & 4 \\ \hline
d & 3 & 2 & 1 & 0 & 3 \\ \hline
x & 2 & 3 & 4 & 3 & 0 \\ \hline
\end{array}
\]
has
\[
[ebc], [ebd], [ecd], [bcd],\; [xeb], [xec], [xdc], [xbc],
\]
but none of $[xed]$, $[exd]$, $[edx]$, $[xbd]$, $[bxd]$, $[bdx]$.

\section{Metric-line equivalences}
\begin{lem}\label{genobs}
If $G$ is a graph such that neither the hypergraph based on $G$ nor the hypergraph based on its complement $\overline{G}$ is metric, then the equivalence relation $\stackrel{G}{\equiv}$ is an obstacle.
\end{lem}
\begin{proof} 
Assuming that $\stackrel{G}{\equiv}$ is not an obstacle, we will prove that at least one of the two hypergraphs based on $G$ and on $\overline{G}$ is metric. This conclusion is immediate when $G$ or $\overline{G}$ is a complete graph, and so we may assume that each of $G$ and $\overline{G}$ has at least one edge. Since $\stackrel{G}{\equiv}$ is not an obstacle, there is a metric space $M$ on a superset $W$ of the vertex set $V$ of $G$ such that
\[
L_M(ab)=L_M(cd) \;\Leftrightarrow\; \{a,b\}\stackrel{G}{\equiv}\{c,d\}
\] 
for every choice of two-point subsets $\{a,b\}$ and $\{c,d\}$ of $V$. 
Let $L$ denote the common value of $L_M(uv)$ with $\{u,v\}$ ranging over the edge set of $G$ 
and let $L'$ denote the common value of $L_M(uv)$ with $\{u,v\}$ ranging over the edge set of $\overline{G}$. Since $L$ and $L'$ are distinct, their symmetric difference $(L-L')\cup (L'-L)$ is nonempty. Switching $G$ and $\overline{G}$ if necessary, we may assume that $L-L'$ is nonempty. Now we are going to prove that the hypergraph based on $G$ is metric. More precisely, with $x$ standing for an arbitrary but fixed element of $L-L'$, we will prove that\\
\phantom{xxx}(i) all three-point subsets of $V$ belong to $\ce_M$,\\
\phantom{xxx}(ii) $x\not\in V$,\\
\phantom{xxx}(iii) if $\{u,v\}$ is an edge of $G$, then $\{x,u,v\}\in\ce_M$,\\
\phantom{xxx}(iv) if $\{u,v\}$ is an edge of $\overline{G}$, then $\{x,u,v\}\not\in\ce_M$.\\
\rule{0pt}{14pt}To prove (i), consider an arbitrary three-point subset $T$ of $V$. Since at least two of the three two-point subsets of $T$ belong to  the same class of $\stackrel{G}{\equiv}$, we may label the elements of $T$ as $u,v,w$ in such a way that 
$\{u,v\}\stackrel{G}{\equiv}\{v,w\}$. This means that $L_M(uv)=L_M(vw)$, and so $w\in L_M(uv)$, and so $\{u,v,w\}\in\ce_M$. 
To prove (ii), we rely on the assumption that each of $G$ and $\overline{G}$ has at least one edge; this assumption along with (i) guarantees that both $L$ and $L'$ contain $V$; now $x\not\in V$ follows from $x\not\in L'$. With (ii) established, (iii) and (iv) follow from $x\in L-L'$.
\end{proof}

\begin{proofof}{Proof of Theorem~\ref{obs}}
Let $n$ be an even integer  greater than four and let $G$ be the cycle $C_n$. 
By Theorem~\ref{evenc}, the hypergraph based on $G$ is not metric; 
by Lemma~\ref{house}, the hypergraph based on $\overline{G}$ is not metric; hence, by Lemma~\ref{genobs}, the
equivalence relation $\stackrel{G}{\equiv}$ is an obstacle. 

To see that $\stackrel{G}{\equiv}$ is a minimal obstacle, consider any proper subset $U$ of $V$.  
The restriction of $\equiv$ on $\binom{U}{2}$ is $\stackrel{F}{\equiv}$, where $F$ is the subgraph of $G$ induced by $U$. Since 
$F$ is an induced subgraph of $P_{n-1}$, Lemma~\ref{oddc} guarantees that there is a metric space $M$ on the ground set $U\cup\{x\}$, where $x\not\in U$, such that all three-point subsets of $U$ belong to $\ce_M$ and such that a two-point subset $\{u,v\}$ of $U$ is an edge of $F$ if and only if $\{x,u,v\}\in \ce_M$. Whenever $\{u,v\}$ is a two-point subset of $U$, we have
\[
L_M(uv)=\begin{cases}
U\cup\{x\}& \text{if $\{u,v\}$ is an edge of $F$},\\
U& \text{otherwise},
\end{cases}
\]
and so $\stackrel{F}{\equiv}$ is not an obstacle.
\end{proofof}

\end{document}